\documentclass[12pt]{article}

\usepackage[utf8]{inputenc}
\usepackage{amsmath,amsfonts,amsthm,amstext,amssymb,amscd,amsbsy}

\usepackage{indentfirst}
\usepackage[dvips]{color}
\usepackage{color}
\usepackage{epsfig}
\usepackage{hyperref}

\usepackage{lscape}
\usepackage{pdflscape}

\usepackage{pstricks, pst-plot}
\usepackage{pictex,dcpic}
\usepackage{tikz}
\usepackage{tikz-cd}

\usepackage{multicol}
\usetikzlibrary{arrows,decorations.markings}
\usetikzlibrary{matrix}
\usetikzlibrary{graphs}
\usetikzlibrary{backgrounds}

\newtheorem{theorem}{Theorem}[section]
\newtheorem{proposition}[theorem]{Proposition}
\newtheorem{lemma}[theorem]{Lemma}
\newtheorem{corollary}[theorem]{Corollary}
\newtheorem{definition}[theorem]{Definition}
\newtheorem{remark}[theorem]{Remark}
\newtheorem{example}[theorem]{Example}

\newcommand{\CC}{\mathbb C}
\newcommand{\RR}{\mathbb R}

\newcommand{\mf}{\mathfrak}
\newcommand{\ad}{\mbox{ad}}
\newcommand{\F}{\mathbb{F}_{\Theta}}
\newcommand{\f}{\mathbb{F}_{\Theta^\prime}}

\newcommand{\p}{\mf{p}}
\newcommand{\m}{\mf{m}}

\begin{document}

%%%%%%%%%%%%%%%%%%%%%%%
% TITLE %
%%%%%%%%%%%%%%%%%%%%%%%
\begin{center} \Large
Geometry of Invariant Almost Semi K\"ahler Submanifolds of Flag Manifolds
\end{center}

\begin{center}
%%%%%%%%%%%%%%%%%%%%%%%
% AUTHOR %
%%%%%%%%%%%%%%%%%%%%%%%
Neiton Pereira da Silva\footnote{Neiton Pereira da Silva at UFU. Email: neiton@ufu.br} 

%%%%%%%%%%%%%%%%%%%%%%%
% FUNDING AGENCY, OPTIONAL %
%%%%%%%%%%%%%%%%%%%%%%%
%Supported by...

%%%%%%%%%%%%%%%%%%%%%%%
% COAUTHORS, OPTIONAL %
%%%%%%%%%%%%%%%%%%%%%%%
%joint work with Coauthor\footnote{Professor at  Y University. Email: y@email.com}
\end{center}

%%%%%%%%%%%%%%%%%%%%%%%
% ABSTRACT %
%%%%%%%%%%%%%%%%%%%%%%%

\begin{abstract}

In this paper we discuss the geometry of homogeneous spaces witch are almost Hermitian submanifolds of flag manifolds. We prove that such spaces are necessarily minimal submanifolds and in the case where these submanifolds are also flag manifolds, they are totally geodesic.

\end{abstract}


%\begin{keywords}

\section{Introduction}

%  It is well known the almost Hermitian manifolds  were classified in 16 classes, by Gray-Hervella in \cite{GH}. Some of them, e.g. Hermitian manifolds, almost K\"ahler manifolds bear sufficient resemblance to K\"ahler manifolds. The several these classes have been studied in flag manifolds, see ....
%  
%  
%  Study geometric properties of almost Hermitian structure that are invariant by the action of a semisimple Lie group $G$, on a flag manifold $G/P$ (here $P$ is a parabolic subgroup of $G$) is very interesting because each definition condition of the 16 classes can be parameterized in terms of roots of the Lie algebra $\mathfrak{g}$ of $G$, see ....
%  Uma propriedade  geométrica que sido estudada é 
% 
%  In \cite{Gray}, Gray 
%
%
%
%falar das classes e artigos relacionados
%
%importância das classes, aplicação em física
%
%importância das flags
%
%teo do Gray
%
%dificuldade de estender o teo de Gray no caso geral
%
%artigo do Watson
%
%caso invariante

%\section{Almost Complex Manifolds}
It is a well known fact that a K\"ahler submanifold is a minimal variety (see \cite{Rham}). This result is also valid for others classes of almost Hermitian structures, \cite{Gray}. In this paper we extend this result for invariant almost semi-K\"ahler submanifolds of flag manifolds. 
  
An almost complex manifold is a manifold $M$ equipped with an tensor field $J$ of type $(1,1)$, called an almost complex structure, such that $J^2=-1$. Such a manifold is orientable and has even dimension, say $2m$. An elementary text on almost complex manifolds can be found in \cite{KN2}.

Recall the Nijenhuis tensor of $J$, defined by:
$$
\dfrac{1}{2}N(X,Y)=[JX,JY]-J[JX,Y]-J[X,JY]-[X,Y], \quad X,Y\in  \mathfrak{X}(M),
$$
a well known theorem of Newlander and Nirenberg states that the almost complex manifold $(M,J)$ is a complex manifold if an only if $J$ is integrable, that is $N=0$. 
An Almost Hermitian manifold  $(M,J,g)$ is an almost complex manifold $(M,J)$ with a $J$-invariant Riemannian metric $g$, i.e., $g(X,Y)=g(JX,JY)$, for any $X,Y \in \mathfrak{X}(M)$. In the sequel we will abbreviate Almost Hermitian manifold by AH-manifold. The K\"ahler form of the almost Hermitian manifold $(M,J,g)$ is the differential form $\Omega(X,Y)=g(X,JY)$, for any $X,Y \in \mathfrak{X}(M)$. 

Following \cite{FIP}, associated to the Levi-Civita connection $\nabla$ on $M$,  we have the following formulas for the covariant derivative, the exterior differentiation and the codifferential of $\Omega$:

$(\nabla_{X}\Omega)(Y,Z)=g(Y,(\nabla_XJ)Z)$,\quad 
$3d\Omega (X,Y,Z)=\sigma(\nabla_X\Omega)(Y,Z)$, and $$(\delta\Omega)(X)=\displaystyle\sum_{i=1}^{m}\left \{\right (\nabla_{E_i}\Omega)(E_i,X)+(\nabla_{JE_i}\Omega)(JE_i,X)\},$$
where $\sigma$ denotes the cyclic sum over $X,Y,Z$ and $\{E_1, \dots,E_m,JE_1, \dots, JE_m\}$ is a local orthonormal J-frame.    

In \cite{GH}, A. Gray and L. M. Hervella classify the AH-manifolds in 16 classes and show certain inclusion relations between them. In the next table we recall some of these classes, with the defining condition.

\begin{tabular}{|c|c|c|}
	\hline
symbol	& name & condition \\
	\hline
K	& K\"ahler & $\nabla J=0$ \\
	%\hline
	AK&Almost K\"ahler  & $d\Omega=0$  \\
	%\hline
NK	& Nearly K\"ahler & $(\nabla_XJ)X =0$\\
%	\hline
QK	& Quasi K\"ahler & $(\nabla_XJ)Y+(\nabla_{JX}J)JY=0$ \\
	%\hline
ASK	& Almost semi-K\"ahler & $\delta\Omega=0$  \\
	%\hline
	SK& Semi-K\"ahler  & $\delta\Omega=0$, $N=0$  \\
	& & (here $N$ is the Nijenhuis tensor)\\
	%\hline
H	& Hermitian &$N=0$  \\
	\hline
\end{tabular}\\

Certains types of almost Hermitian manifolds ( e.g., H, AK) have been studied by several authors with the aim of generalizing the K\"ahler geometry. 

A $2r$-dimension submanifold $M^\prime$ of $M$ is called holomorphic (or invariant or almost-complex) submanifold of $M$ if for any $p\in M^\prime$ the tangent space $T_pM^\prime$ is $J$-invariant, i.e. $J(T_pM^\prime)=T_pM^\prime$.

In \cite{Gray}, A. Gray obtained:\textit{ If $M$ belongs to one of the classes QK, AK, NK, H, K, then any holomorphic submanifold $M^\prime$ of $M$ belongs to the same class. Moreover, $M^\prime$ is necessarily mininal, except in the Hermitian case.}

Note that the classes  ASK and SK  are excluded from the above Gray´s result. According \cite{FIP}, in general, it is not known if holomorphic submanifolds of almost semi-K\"ahler manifolds are minimal. However, if $M^\prime$ is a holomorphic submanifold of an ASK-manifold with codimension 2, then $M^\prime$ is minimal (see \cite{Watson}). 

In this paper we show that Gray's result can be extended for the classes ASK and SK when $M=G/P$ is a flag manifold (i.e. $M$ is a homogeneous space $G/P$ where $G$ is semi-simple Lie group and $P$ is a parabolic subgroup) and $G/P$ is endowed with a almost Hermitian structure $(J,g)$ invariant by the action of $G$. More precisely, we prove that if $(G/P, J,g)$ is an invariant almost Hermitian flag manifold (witch belong to the class ASK or SK). Then any  invariant holomorphic submanifold $M^\prime$  of  $G/P$ belongs to the same class. Moreover, $M^\prime$ is necessarily minimal. See Theorem (\ref{ask}) and Corollary (\ref{sk}). We use the linear operators defined in \cite{FIP} to study the minimally of $M^\prime$. 

In the case that an invariant holomorphic submanifold $M^\prime$  of  $G/P$ is also a flag manifold, we obtain an very interesting result witch states that $M^\prime$ is a totally geodesic submanifold with respect to any invariant metric $g$ on $G/P$, see Theorem \ref{tg}. Here we used properties of Lie theory to obtain the result. 

Others authors have studied (and have been studying) geometric properties of some classes of almost Hermitian structures on flag manifolds, see for example \cite{SM N}, \cite{SM R} and \cite{LN}. The main idea here is that the definition conditions of the 16 classes can be parameterized in terms of roots of the Lie algebra $\mathfrak{g}$ of $G$. Thus investigating  geometric and others properties reduces to a study in terms of Lie theory. 

In section 2, we recall the linear operators and others results used to obtain the mains theorems of this paper. In section 3, we give the description of flag manifold as a even dimension real homogeneous space by means of Lie theory. In section 4, we set our notation for invariant almost Hermitian structures on flag manifolds recalling the notion of invariant almost complex structure and Riemannian metric on a flag manifold. In section 5, we compute explicitly the Riemannian connection of invariant Riemannian metric on the Weyl basis. In section 6, we that the covariant derivative of the K\"ahler form is zero when computed in a invariant global orthonormal J-frame. As an consequence we obtain that all flag manifolds belong to the class ASK. In section 7, we study homogeneous holomorphic submanifolds of flag manifolds. One of the main result of this paper is Theorem \ref{ask}. In section 8, we study the case where the homogeneous holomorphic submanifolds is also a flag manifold and obtain the other main result, Theorem \ref{tg}.

\section{Almost Complex submanifolds}

Now let $(M,J,g)$ be an AH-manifold and $M^\prime$ a holomorphic submanifold.  The almost complex structure induced on $M^\prime$ and the almost Hermitian Riemannian metric induced on $M^\prime$ will be also denoted by $J$ and $g$, respectively. Note that $J(T_p M^{\prime\perp})=T_pM^{\prime\perp}$, because if $Z\in T_pM^{\prime\perp} $ then $g(JZ,X)=-g(Z,JX)=0$, for any $X \in T_pM^\prime$, since $J(T_pM^\prime)=T_pM^\prime$.

Recall the Gauss formula:

$$
\nabla_XY=\nabla^\prime_XY+\alpha(X,Y), \quad\quad X,Y\in \mathfrak{X}(M^\prime)
$$
where $\nabla^\prime$ and $\alpha$  denote the Riemannian connection and the second fundamental form of $M^\prime$, respectively.

%\begin{remark}

%\end{remark}
%Now, for any 

An important type of holomorphic submanifold are the minimal one, in particular those totally geodesics. Next we recall these notions, see \cite{O'Neil}, p. 101,or \cite{Besse} p.38, for details.

\begin{definition}
	Let $(M^\prime, g,J)$ be an almost complex submanifold of $(M, g,J)$. The mean curvature of $M^\prime$ at $x\in M^\prime$ is the normal vector
	$$
	H_x=\dfrac{1}{2r}\displaystyle\sum_{i=1}^{r}(\alpha(E_i,E_i)+\alpha(JE_i,JE_i))
	$$ 
where $\{E_1, \dots,E_r,JE_1, \dots, JE_r\}$ is a local orthonormal J-frame on $TM^\prime$,
\begin{enumerate}
	\item[a)]$(M^\prime, g,J)$ is minimal if $H\equiv 0$; 
	\item[b)] $(M^\prime, g,J)$ is totally geodesic if $\alpha\equiv 0$.
\end{enumerate}
 
\end{definition}

Next, we denote by $\bar{\mathfrak{X}}(M)$ the Lie algebra of vector fields tangent to $M$ along $M^\prime$. Thus,
$$
\mathfrak{X}(M^\prime)\subset \bar{\mathfrak{X}}(M)\subset \mathfrak{X}(M),
$$
and 
\begin{equation}\label{fields perp}
\bar{\mathfrak{X}}(M)=\mathfrak{X}(M^\prime)\oplus \mathfrak{X}(M^\prime)^\perp.
\end{equation}

\begin{definition}
	Let $M^\prime$ be a holomorphic submanifold of an AH-manifold $(M,J,g)$. The partial (or tangent) coderivative of $\Omega$ and the normal-coderivative of $\Omega$, with respect to $M^\prime$ are the linear operators $\bar{\delta}\Omega$ and $\bar{\bar{\delta}}\Omega$, respectively, given by  :  
	
	$$
	(\bar{\delta}\Omega)(X)=\displaystyle\sum_{i=1}^{r}\left \{\right (\nabla_{E_i}\Omega)(E_i,X)+(\nabla_{JE_i}\Omega)(JE_i,X)\},
	$$
	$$
	(\bar{\bar{\delta}}\Omega)(X)=\displaystyle\sum_{i=1}^{m-r}\left \{\right (\nabla_{F_i}\Omega)(F_i,X)+(\nabla_{JF_i}\Omega)(JF_i,X)\},
	$$
	where $X\in \bar{\mathfrak{X}}(M) $;  $\{E_1, \dots,E_r,JE_1, \dots, JE_r\}$ is a local orthonormal J-frame on $TM^\prime$, $\dim M^\prime=2r$ and $\{F_1, \dots,F_{m-r},JF_1, \dots, JF_{m-r}\}$ is a local orthonormal J-frame on $\left( TM^\prime\right) ^\perp$.   
\end{definition}

An easy computation shows that by using the Gauss formula,
\begin{equation}\label{gauss}
	\nabla_{X}J=\nabla^\prime_{X}J+\alpha_{X}J, \quad X\in\mathfrak{X}(M^\prime)
\end{equation}
where $(\alpha_{X}J)Y=\alpha(X,JY)-J\alpha(X,Y)$. Using equation (\ref{gauss}), one proves that for any $Z\in \mathfrak{X}(M^\prime)^\perp$:
\begin{equation}\label{minimal}
(\bar{\delta}\Omega)(Z)=2rg(JH,Z)
\end{equation}
where $H$ is the mean curvature vector field of the submanifold $M^\prime$ (see \cite{FIP} p. 63 or \cite{Gray}).
%Indeed,  from (\ref{gauss})  one has
%\begin{eqnarray*}
%(\bar{\delta}\Omega)(Z)&=&\displaystyle\sum_{i=1}^{r}\left \{\right (\nabla_{E_i}\Omega)(E_i,Z)+(\nabla_{JE_i}\Omega)(JE_i,Z)\}\\
%&=&\displaystyle\sum_{i=1}^{r}\left \{ g(E_i,(\nabla_{E_i}J)Z)+g(JE_i,(\nabla_{JE_i}J)Z)\right\}\\
%&=&\displaystyle\sum_{i=1}^{r}\left \{ g(E_i,(\nabla^\prime_{E_i}J+\alpha_{E_i}J)Z)+g(JE_i,(\nabla^\prime_{JE_i}J+\alpha_{JE_i}J)Z)\right\}\\
%&=&\displaystyle\sum_{i=1}^{r}\left \{ g(E_i,(\alpha_{E_i}J)Z)+g(JE_i,(\alpha_{JE_i}J)Z)\right\}\\
%&=&\displaystyle\sum_{i=1}^{r}\left \{ g(J(\alpha(E_i,E_i))+J(\alpha(E_i,JE_i)),Z)+g(J(\alpha(JE_i,JE_i))+J(\alpha(JE_i,-E_i)),Z)\right\}\\
%&=& g(\displaystyle\sum_{i=1}^{r}\left \{ J(\alpha(E_i,E_i))+J(\alpha(JE_i,JE_i))\right\},Z)\\
%&=&2rg(JH,Z).\qed
%\end{eqnarray*}
In particular, a almost complex submanifold $M^\prime$ of $(M,J,g)$ is minimal if and only if $(\bar{\delta}\Omega)=0$ on $\mathfrak{X}(M^\prime)^\perp$.

\begin{proposition}\cite{FIP}
Let $M^\prime$ be a holomorphic submanifold of an AH-manifold $(M,J,g)$. Then, for any $X\in \mathfrak{X}(M^\prime) $, one has:
\begin{equation}\label{same class}
({\delta}\Omega)(X)=({\delta}^\prime\Omega)(X)+ (\bar{\bar{\delta}}\Omega)(X),
\end{equation}
 $\delta^\prime$ denoting the codifferential operator on $M^\prime$.
\end{proposition}

\section{Flag manifolds}

In this section we set up our notation and present the standard
theory of partial (or generalized) flag manifolds associated with semisimple
Lie algebras (see, for example, \cite{LN} or \cite{Arv}  for a similar description).

Let $\mathfrak{g}$ be a finite-dimensional semisimple complex Lie
algebra and $G$ a Lie group with Lie algebra $\mathfrak{g}$. Consider a Cartan subalgebra $\mathfrak{h}$ of $\mathfrak{g}$. We denote by $R$ the system of
roots of $(\mathfrak{g},\mathfrak{h})$. A root $\alpha\in R$ is a linear
functional on $\mathfrak{g}$. It determines uniquely an element
$H_{\alpha}\in\mathfrak{h}$ by the Riesz representation $\alpha(X)=B(
X,H_\alpha)$, $X\in\mathfrak{g}$, with respect to the Killing form
$B(\cdot,\cdot)$ of $\mf{g}$. The Lie algebra $\mathfrak{g}$ has the following
decomposition
\[
\mathfrak{g}=\mathfrak{h}\oplus\sum_{\alpha\in R}\mathfrak{g}_{\alpha}
\]
where $\mathfrak{g}_{\alpha}$ is the one-dimensional root space corresponding to $\alpha$. The eigenvectors $E_{\alpha}\in\mathfrak{g}_{\alpha}$ satisfy the following equation

\begin{equation}
	\left[  E_{\alpha},E_{-\alpha}\right]  =B\left(
	E_{\alpha},E_{-\alpha}\right) H_{\alpha}. \label{igualdade}
\end{equation}

We fix a system $\Sigma$ of simple roots of $R$ and denote by
$R^+$ and $R^{-}$ the corresponding set of positive and
negative roots, respectively. Let $\Theta\subset\Sigma$ be a
subset, define
\begin{align}
	R(\Theta):=\langle\Theta\rangle\cap R \quad \text{and}\quad
	R(\Theta)^\pm:=\langle\Theta\rangle\cap R^{\pm}.\nonumber
\end{align}
We denote by $R_{\Theta}:=R\setminus R(\Theta)$ the complementary set of roots. In general, $R_{\Theta}$ is not a root system, see an example in \cite{LN}. 

%Recall that a root system $R$ is irreducible if and only if $R$ (or, equivalently, $\Sigma$) cannot be partitioned into two proper, orthogonal subsets (see \cite{H}). Equivalently, $R$ is irreducible iff the Dynkin diagram of $\Sigma$ is connected. Besides that it is well known if $\mf{g}$ is a simple Lie algebra then its root system $R$ is irreducible. 

The parabolic subalgebra $\mathfrak{p}_{\Theta}$, associated to $\Theta$, is defined by
\[
\p_{\Theta}:=\mf{h}\oplus\sum_{\alpha\in R^+}\mf{g}_{\alpha}\oplus\sum_{\alpha\in R(\Theta)^{-}}\mf{g}_{\alpha}.
\]
Note that $\p_{\Theta}$ contains the Borel subalgebra
$\mf{b}^+=\mf{h}\oplus\sum\limits_{\alpha\in R^+}\mf{g}_{\alpha}$.

The partial flag manifold determined by the choice of $\Theta\subset \Sigma$ is the
homogeneous space $\F=G/P_{\Theta}$, where $P_{\Theta}$ is the normalizer of
$\mf{p}_{\Theta}$ in $G$. In the special case $\Theta=\emptyset$, we obtain the maximal flag manifold $\mathbb{F}=G/B$ associated with $\mf{g}$, where
$B$ is the normalizer of the Borel subalgebra.

Now we will see the construction of any flag manifold as the
quotient $U/K_{\Theta}$ of a  compact semisimple Lie group
$U\subset G$ modulo the centralizer $K_{\Theta}$ of a torus in
$U$. We fix once and for all a Weyl basis of $\mf{g}$ which amounts
to giving $X_\alpha\in\mf{g_{\alpha}}$, $H_{\alpha}\in\mf{h}$ for $\alpha\in R$, with the standard
properties:
\begin{equation}\label{base weyl}
	\begin{tabular}
		[c]{lll}
		$B( X_{\alpha},X_{\beta}) =\left\{
		\begin{array}
			[c]{cc}%
			1, & \alpha+\beta=0,\\
			0, & \text{otherwise};
		\end{array}
		\right. $ &\hspace{-.5cm} & $\left[  X_{\alpha},X_{\beta}\right]  =\left\{
		\begin{array}
			[c]{cc}%
			H_{\alpha}\in\mathfrak{h}, & \alpha+\beta=0,\\
			n_{\alpha,\beta}X_{\alpha+\beta}, & \alpha+\beta\in R,\\
			0, & \text{otherwise.}%
		\end{array}
		\right.  $
	\end{tabular}
\end{equation}
The real constants $n_{\alpha,\beta}$ are non zero if and only if
$\alpha+\beta\in R$ and satisfy
$$
\left\{
\begin{array}
	[c]{cc}%
	n_{\alpha,\beta}=-n_{-\alpha,-\beta}=-n_{\beta,\alpha}& \\
	\hspace*{-1cm} n_{\alpha,\beta}=n_{\beta,\gamma}=n_{\gamma,\alpha},&\mbox{if}\quad\alpha+\beta+\gamma=0.
\end{array}
\right.
$$

Consider the following two-dimensional real spaces
%\begin{align}
$\mf{u}_{\alpha}=\text{span}_\mathbb{R}\{A_{\alpha},iS_{\alpha}\}$, 
%,\hspace{0.3cm}
with $A_{\alpha}=X_{\alpha}-X_{-\alpha}$
%\hspace{0.3cm}
and $S_{\alpha}=X_{\alpha}+X_{-\alpha}$, 
$\alpha\in R^{+}$.
%\end{align}
Then the real Lie algebra
$\mf{u}=i\mf{h}_\mathbb{R}\oplus\sum\mf{u}_{\alpha}$,  $\alpha\in R^{+},$
is a compact real form of $\mf{g}$, where $\mf{h}_\mathbb{R}$ denotes the real vector space spanned by $\{H_\alpha: \,\alpha\in R\}$.

Let $U=\exp\mf{u}$ be the compact real form of $G$ corresponding
to $\mf{u}$. By the restriction of the action of $G$ on $\F$, we
can see that $U$ acts transitively on $\F$, then $\F=U/K_{\Theta}$,
where $K_{\Theta}=P_{\Theta}\cap U$. The Lie algebra $\mf{k}_{\Theta}$ of $K_\Theta$ is the set of fixed points of the conjugation $\tau\colon X_{\alpha}\mapsto-X_{-\alpha}$ of $\mf{g}$ restricted to $\p_{\Theta}$
$$
\mf{k}_{\Theta}=\mf{u}\cap\p_{\Theta}=i\mf{h}_\mathbb{R}\oplus\sum_{\alpha\in R(\Theta)^{+}}\mf{u}_{\alpha}.
$$

The tangent space of $\F=U/K_{\Theta}$ at the origin $o=eK_{\Theta}$ can be identified with the orthogonal complement (with respect to the Killing form) of $\mf{k}_{\Theta}$ in $\mf{u}$
$$
T_o\F=\m=\sum\limits_{\alpha\in R_{\Theta}^{+}}\mf{u}_{\alpha},
$$
with $R_{\Theta}^{+}=R_{\Theta}\cap R^+$. Thus we have $\mf{u}=\mf{k}_{\Theta}\oplus\mf{m}$. The complexification of $\m$ is given by:
$$
\m^\CC=\m \otimes  \CC= \sum\limits_{\alpha\in R_{\Theta}}\CC X_{\alpha}.
$$ 

\section{Invariant almost Hermitian structures on flag manifolds}

\subsection{Invariant almost complex structures}
In this section we describe almost complex structure on $\F$ and Riemannian metric, both invariant by the action of the Lie group $U$, similiar description can be found in \cite{LN}, \cite{SM R}, or \cite{SM N} for maximal flag manifolds. 

An \emph{almost complex structure} on $\F=U/K_{\Theta}$ is a tensor field of type $(1,1)$ that corresponds each $x\in \F$ to a linear endomorphism  $J_{x}\colon T_{x}\F\rightarrow T_{x}\F$ which satisfies $J_{x}^{2}=-\operatorname{Id}$. The almost complex structure on $\F$ is invariant (or $U$-invariant) if 
\[
d \mathit{u}_{x}\circ J_{x}=J_{\mathit{u}x}\circ d\mathit{u}_{x}
\]
for all $\mathit{u}\in U$. An \textit{invariant almost complex structure} (\textit{iacs} from now on) is determined by a linear endomorphism $J\colon\mathfrak{m\longrightarrow} \mathfrak{m}$, which satisfies $J^{2}=-\operatorname{Id}$ and commutes with the adjoint action of $K_\Theta$ on $\mathfrak{u}$, that is,
\[
\operatorname{Ad}(k) J=J \operatorname{Ad}(k), \quad \text{for all}\quad  k\in K_\Theta,
\]
or, equivalently,
\[
\operatorname{ad}(L)J=J\operatorname{ad}(L), \quad \text{ for all} \quad L\in \mathfrak{k}_\Theta.
\]

We will use the same letter $J$ to denote its extension to the complexification $\mathfrak{m}^{\mathbb{C}}$. Since $J^{2}=-\operatorname{Id}$, its eigenvalues are $i$ and $-i$ and the corresponding eigenspaces are denoted by
\[
T^{(1,0)}_{o}\F=\left\{X\in T_{o}\mathbb{F}^{\mathbb{C}}:
JX=iX\right\}
\]
and
\[
T^{(0,1)}_{o}\F=\left\{X\in T_{o}\mathbb{F}^{\mathbb{C}}:
JX=-iX\right\}.
\]

Thus
\[
\mathfrak{m}^{\mathbb{C}}=T_{o}\F^{\mathbb{C}}=T^{(1,0)}_{o}\F\oplus
T^{(0,1)}_{o}\F.
\]
The eigenvectors with eigenvalue $+i$ (resp. $-i$) are called  \textit{of type (1,0)} (resp. \textit{of type (0,1)}). 

The invariance of $J$ entails that $J(\mf{g}_\alpha)=\mf{g}_\alpha$ for all $\alpha\in R_\Theta$. The eigenvalues of $J$ are $\pm i$ and the eigenvectors in $\mf{m}^\CC$ are $X_\alpha$; $\alpha\in R_\Theta$. Hence 
$$
J(X_\alpha)= i\varepsilon_\alpha X\alpha
$$
with $\varepsilon_\alpha=\pm 1$ satisfying $\varepsilon_{-\alpha}=-\varepsilon_\alpha$.

Thus an iacs on $\F$ is completely prescribed by a set of signs $\{\varepsilon_{\alpha}\}_{\alpha\in R_\Theta}$ with $\varepsilon_{-\alpha}=-\varepsilon_\alpha$. %In the sequel we abuse notation and say that an invariant almost complex structure on $\F$ is $J=\{\varepsilon_\alpha\}$.
Since $\F$ is a homogeneous space of a complex Lie group it has a natural structure
of a complex manifold. The associated integrable almost complex structure $J_\CC$ is
given by $\varepsilon_\alpha=1$ if $\alpha>0$. The conjugate structure $-J_\CC$ is also integrable. A complete description of complex structures of flag manifolds can be found in \cite{Alek-Arv}.  %The next result classifies the invariant complex structures on flag manifolds.

\subsection{Riemannian invariant metrics}

Now we will describe Riemannian metrics $g(\cdot,\cdot)$ on $\F$ which is $J$-invariant, i.e. $g(JX,JY)=g(X,Y)$ and also invariant by the action of $U$ (or $G$).  It is well known a Riemannian invariant metric on $\F$ is completely determined by a real inner product $g\left(\cdot,\cdot\right)$ on $\mathfrak{m}=T_{o}\F$ which is invariant by the adjoint action of $\mf{k}_{\Theta}$. Besides that any real inner product $\ad(\mf{k}_{\Theta})$-invariant on $\mf{m}$ has the form

\begin{equation}
	g\left(\cdot,\cdot\right)=-\displaystyle\sum_{\alpha\in R_{\Theta}^+}\lambda_\alpha B\left(\cdot,\cdot\right)|_{\mathfrak{u}_{\alpha}\times\mathfrak{u}_{\alpha}}
\end{equation}
where the Killing form $B(\cdot,\cdot)$ is negative defined on $\mf{m}\times \mf{m}$ and $\lambda_\alpha>0$. %We will denote a Riemannian invariant metric on $\F$  by $g\left(\cdot,\cdot\right)$.

Since $\F$ is a reductive homogeneous space,  i.e. $\operatorname{ad}(\mf{k}_{\Theta})\mf{m}\subset\mf{m}$, we can decompose $\mf{m}$ into a sum of irreducible $\ad (\mf{k}_{\Theta})$ submodules $\mf{m}_i$ of the module $\mf{m}$:
\[
\m= \m_1\oplus\cdots\oplus\m_s.
\]

%Now we describe the invariant metrics on flag manifolds.
Thus a Riemannian invariant metric on $\F$ is completely determined by a real inner product $g\left(\cdot,\cdot\right)$ on $\mathfrak{m}=T_{o}\F$ which is invariant by the adjoint action of $\mf{k}_{\Theta}$. From Schur`s Lemma any real inner product $\ad(\mf{k}_{\Theta})$-invariant on $\mf{m}$ is given by

\begin{equation}\label{inner product}
	g\left(\cdot,\cdot\right)=-\lambda_1 B\left(\cdot,\cdot\right)|_{\mathfrak{m}_{1}\times\mathfrak{m}_{1}}-\cdots
	-\lambda_s B\left(\cdot,\cdot\right)|_{\mathfrak{m}_{s}\times\mathfrak{m}_{s}}.
\end{equation}\\
Here $\lambda_i=\lambda_\alpha$ if $\mf{u}_\alpha\subset \mf{m}_i$, hence $\lambda_\alpha=\lambda_\beta$ if $\mf{u}_\alpha$ and $\mf{u}_\alpha$ belongs to the same irreducible $\ad (\mf{k}_{\Theta})$ submodules $\mf{m}_i$ of the module $\mf{m}$.
%where ${\mathfrak{m}}_{i}=\mf{m}_{\xi_i}$ and $\lambda_i=\lambda_{\xi_i}>0$ with $\xi_{i}\in R^{+}_{\mathfrak{t}}$, for $i=1,\ldots,s$. 
%So any invariant Riemannian metric on $\F$ is determined by $|R_{\mf{t}}^+|$ positive parameters.
We will call an inner product defined by equation (\ref{inner product}) as a Riemannian invariant metric on $\F$. Any such inner product  extends naturally to a symmetric bilinear form on the complexification $\mathfrak{m}^{\mathbb{C}}$ of $\mathfrak{m}$, we will denote this extension by the same notation $g\left(\cdot,\cdot\right)$.

% As usual the eigenvectors associated to þi
%are said to be of type ð1; 0Þ while the i-eigenvectors are of type ð0; 1Þ: Thus the
%ð1; 0Þ vectors are multiples of Xa; ea ¼ þ1; and the ð0; 1Þ multiples of Xa; ea ¼ 1:

%The iacs was classified in terms of roots of $R$ for maximal flags manifolds in \cite{SM N} and for partial ones in \cite{SM R}. In the next result we classify iacs in terms of $\mf{t}$-roots.

%The next result is known for maximal flag manifolds (see \cite{SM N}). A similar result for partial flag manifolds can be found in \cite{SM R}. Our work here is to classify iacs in terms of $\mathfrak{t}$-roots. 
%\begin{remark}
	Using the Weyl basis chosen it is easy to see that any iacs $J$ is compatible with any invariant metric $g$ on $\F$, i.e. $g(JX,JY)=g(X,Y)$, for any $X,Y\in \mf{m}$. The pair $(g,J)$ is called an \textit{invariant almost Hermitian structure} on $\F$.
%\end{remark}

%An invariant metric of the form $g(\cdot, \cdot)=-\lambda B\left(\cdot,\cdot\right)|_{\mathfrak{m}}$ is called a normal metric on $\F$.

\begin{remark}\label{J-orthonomal frame}
	By an easy computation, one see that an (invariant) $J$ frame on $\F$ is given by 
$$
	\left\lbrace A_\alpha/\sqrt{2\lambda_\alpha}, J(A_\alpha/\sqrt{2\lambda_\alpha})=i\varepsilon_\alpha S_\alpha/\sqrt{2\lambda_\alpha}; \quad \alpha \in R_{\Theta}^+ \right\rbrace. 
$$
\end{remark}

\section{The invariant Riemannian connection on $\F$}

In this section we compute explicitly the invariant Levi-Civita connection on $\F$.

Recall the isomorphism $\m\rightarrow T_o\F$ which corresponds each $X\in\mf m$ to the  invariant vector field $X^\star\in \mathfrak{X}(\F)$:
$$
X^\star_o =\dfrac{d}{dt}(exp(tX))\cdot o\left| _{t=0}\right..  
$$

For each $X\in\mf m$ the vector field $X^\star$ is a Killing vector field, since the Riemannian metric $g(\cdot, \cdot)$ is invariant. This is equivalent to  
\begin{equation}\label{Xg}
	X^\star g(Y,Z)=g([X^\star,Y], Z)+g(Y,[X^\star,Z])
\end{equation}
and 
\begin{equation}\label{antisimet}
	g(\nabla_YX^\star,Z)+g(\nabla_ZX^\star,Y)=0
\end{equation}

for all $Y,Z \in \mathfrak{X}(\F)$, (see \cite{O'Neil}, pp.251).

Now recall that if $X,Y \in \mf m$ then $\left[X^\star,Y^\star \right]=-\left[X,Y \right]^\star $. Then using (\ref{Xg}) and (\ref{antisimet}), it can be proved that the  Riemannian connection is given by (\cite{Besse}, p. 183 or \cite{Arv book}, p. 79) :

$$(\nabla_{X^{\star}} Y^\star)_o=-\dfrac{1}{2}[X,Y]_{\mf m}+U(X,Y), \quad X,Y \in \mf m$$

where $U\colon \mf m\times\mf m\longrightarrow \mf m$ is the symmetric bilinear map defined by 

\begin{equation}\label{U}
2g(U(X,Y),Z)= g([Z,X]_{\mf m},Y)+g(X, [Z,Y]_{\mf m}).
\end{equation}

We also denote $U(X,Y)$ its complexification $U\colon \mf{m}^\CC\times\mf{m}^\CC\longrightarrow \mf{m}^\CC$. On the basis $\{X_\alpha\}$, $\alpha \in R_{\Theta}$,  $U(X_\alpha, X_\beta)=\sum a_\gamma X_\gamma$ then
$$
g(U(X_\alpha, X_\beta),X_{-\gamma})=a_\gamma\lambda_\gamma.
$$
Now 
$
[X_\alpha, X_\beta]_{\mf m}\neq 0 \iff \alpha+\beta \in R_{\Theta}=R\setminus R(\Theta)
$ and in this case $[X_\alpha, X_\beta]_{\mf m}=n_{\alpha,\beta}X_{\alpha+\beta}$.

Thus, 

\begin{eqnarray*}
2g(U(X_\alpha, X_\beta),X_-\gamma)&=&g([X_-\gamma,X_\alpha]_{\mf m},X_\beta)+g(X_\alpha,[X_-\gamma,X_\beta]_{\mf m})\\&=&\left\{
\begin{array}
	[c]{cc}%
	n_{\alpha,\beta}(\lambda_\beta-\lambda_\alpha), & \mbox{if}\quad \gamma=\alpha+\beta \\ \\
	\hspace*{-1cm} 0,& \mbox{otherwise.}
\end{array}
\right.
\end{eqnarray*}
Then
$$
U(X_\alpha, X_\beta)=\left\{\begin{array}
	[c]{cc}%
	\dfrac{n_{\alpha,\beta}}{2\lambda_{\alpha+\beta}}(\lambda_\beta-\lambda_\alpha)X_{\alpha+\beta}, & \mbox{if}\quad \alpha+\beta\in R_{\Theta} \\ \\
	\hspace*{-1cm} 0,& \mbox{otherwise.}
\end{array}
\right.
$$

Thus we obtain the next result.
\begin{proposition}
The Levi-Civita connection for invariant vector fields on $\F$ is given by
\begin{eqnarray*}
	\left( \nabla_{X^\star_\alpha} X^\star_\beta\right)_o &=&-\dfrac{1}{2}[X_\alpha,X_\beta]_{\mf m}+U(X_\alpha,X_\beta)\\\\&=&
	\left\{\begin{array}
		[c]{cc}%
		\dfrac{n_{\beta,\alpha}}{2\lambda_{\alpha+\beta}}(\lambda_{\alpha+\beta}+\lambda_\alpha-\lambda_\beta)X_{\alpha+\beta}, & \mbox{if}\quad \alpha+\beta\in R_{\Theta} \\ \\
		\hspace*{-1cm} 0,& \mbox{otherwise.}
	\end{array}
	\right.
\end{eqnarray*}
\end{proposition}
%Note that  $\left( \nabla_{X^\star_\alpha} X^\star_\beta\right)_o-\left( \nabla_{X^\star_\beta} X^\star_\alpha\right)_o=-\left[ X_\alpha,X_\beta\right] =- \left[  X_\alpha ,  X_\beta\right]^\star_o=\left[  X^\star_\alpha ,  X^\star_\beta\right]_o$ and $\left( X_\alpha^\star g(X_\beta^\star,X_{-(\alpha+\beta)}^\star)\right) _o=g(\nabla_{X^\star_\alpha} X^\star_\beta,X_{-(\alpha+\beta)}^\star)_o+g(X_\beta^\star,\nabla_{X^\star_\alpha} X^\star_{-(\alpha+\beta)})_o$, according to equation (\ref{Xg}).

For later using, we denote 
$$
r_{\alpha, \beta}=\dfrac{n_{\beta,\alpha}}{2\lambda_{\alpha+\beta}}(\lambda_{\alpha+\beta}+\lambda_\alpha-\lambda_\beta)
$$
then
\begin{equation}\label{connection}
\left( \nabla_{X^\star_\alpha} X^\star_\beta\right)_o=r_{\alpha, \beta} X_{\alpha+\beta}.
\end{equation}
A direct computation shows that
\begin{eqnarray*}
	r_{\alpha, \beta}&=&-r_{-\alpha, -\beta};\\ \\ r_{\alpha, \beta}+r_{\beta,\alpha}&=&\dfrac{n_{\beta,\alpha}}{\lambda_{\alpha+\beta}}(\lambda_\alpha-\lambda_\beta).\\
%	r_{\alpha, \beta}=r_{\beta,\gamma}&=& r_{\gamma,\alpha}\quad \text{if}\quad \alpha+\beta+\gamma=\dfrac{n_{\beta,\alpha}}{2}\left[3+\dfrac{\lambda_\beta(\lambda_\alpha^2-\lambda_\gamma^2)+\lambda_\alpha(\lambda_\gamma^2-\lambda_\beta^2)+\lambda_\gamma(\lambda_\beta^2-\lambda_\alpha^2)}{\lambda_\alpha\lambda_\beta\lambda_\gamma} \right]. 
\end{eqnarray*}

\begin{lemma} \label{connections AS}
	Let $\alpha, \beta$  be roots in $ R_{\Theta} $ such that $\alpha+\beta$ and $\alpha-\beta$ are in $ R_{\Theta} $. Then 
	\begin{enumerate}
		\item[i)] $\nabla_{S_\beta}A_\alpha=r_{\beta,\alpha}S_{\alpha+\beta}-r_{\beta,-\alpha}S_{\beta-\alpha}$
		\item[ii)] $\nabla_{S_\beta}S_\alpha=r_{\beta,\alpha}A_{\beta+\alpha}+r_{\beta,-\alpha}A_{\beta-\alpha}$
		\item[iii)] $\nabla_{A_\beta}A_\alpha=r_{\beta,\alpha}S_{\beta+\alpha}-r_{\beta,-\alpha}A_{\beta-\alpha}$
		\item[iv)] $\nabla_{A_\beta}S_\alpha=r_{\beta,\alpha}A_{\beta+\alpha}+r_{\beta,-\alpha}S_{\beta-\alpha}$
	\end{enumerate} 
	where $r_{\beta,\alpha}=\dfrac{n_{\alpha,\beta}}{2\lambda_{\alpha+\beta}}(\lambda_{\alpha+\beta}+\lambda_{\beta}-\lambda_{\alpha})$.
\end{lemma}

\begin{proof}
All items follows by a direct computations, using    $A_{\alpha}=X_{\alpha}-X_{-\alpha}$, $S_{\alpha}=X_{\alpha}+X_{-\alpha}$ and equation (\ref{connection}).
\end{proof}

%Its is well known that $\F$ admits a K\"ahler-Einstein metric $g$, with respect to the natural complex structure $J_{\mathbb{C}}$.  particular, the K\"ahler connection, i.e. the Riemannian connection with the K\"ahler metric, is given by:
%$$
%\left( \nabla_{X^\star_\alpha} X^\star_\beta\right)_o= \left\{\begin{array}
%	[c]{cc}%
%	n_{\beta,\alpha}\dfrac{\lambda_{\alpha}}{\lambda_{\alpha+\beta}}X_{\alpha+\beta}, & \mbox{if}\quad \alpha+\beta\in R_{\Theta} \\ \\
%	\hspace*{-1cm} 0,& \mbox{otherwise.}
%\end{array}
%\right.
%$$

\section{Almost semi K\"ahler Flag manifolds}

In this section we show the covariant derivative of the K\"ahler form $\Omega$ has a special property, as an immediate consequence one obtains that $\F$ always belong to class ASK for any invariant almost Hermitian structure $(g,J)$.

Consider the codifferential of $\Omega$ on a almost Hermitian  Riemannian manifold $(M, g, J)$:
$$
(\delta\Omega)(X)=\sum_{i=1}^{m}((\nabla_{E_i}\Omega)(E_i,X)+(\nabla_{JE_i}\Omega)(JE_i,X)),
$$
where $\{E_1,\dots, E_m,JE_1,\dots,JE_m\}$ is a local orthonormal $J$-frame. Recall that $(M, g, J)$ is called Almost Semi K\"ahler (ASK) if $\delta\Omega=0$.

In the case of flag manifolds, given a iacs $J$ on $\F$, from the invariance of the vector fields, from remark \ref{J-orthonomal frame}, a global orthonormal $J$-frame at the origin is given by
\begin{equation}\label{J-frame}
	\left\lbrace V_\alpha= A_\alpha/\sqrt{2\lambda_\alpha}, J(V_\alpha)=J(A_\alpha/\sqrt{2\lambda_\alpha})=i\varepsilon_\alpha S_\alpha/\sqrt{2\lambda_\alpha}; \quad \alpha \in R_{\Theta}^+ \right\rbrace .
\end{equation}

\begin{theorem}\label{codifferential}
	For any invariant almost Hermitian structure on $\F$, one has
	$$ (\nabla_{V_\beta}\Omega)(V_\beta,X)=(\nabla_{JV_\beta}\Omega)(JV_\beta,X)=0
	$$
	where $X\in \mf{m}$.
\end{theorem}

\begin{proof}

Using the linearity of the invariant metric, the iacs and the Riemmanian connection, it is easy to see that if   
$$
X=\sum_{\alpha \in R_{\Theta}^+} a_\alpha V_\alpha +\sum_{\alpha \in R_{\Theta}^+} b_\alpha JV_\alpha
$$ then
$$
(\nabla_{V_\beta}\Omega)(V_\beta,X)=\sum_{\alpha \in R_{\Theta}^+} a_\alpha (\nabla_{V_\beta}\Omega)(V_\beta,V_\alpha)+\sum_{\alpha \in R_{\Theta}^+} b_\alpha (\nabla_{V_\beta}\Omega)(V_\beta,JV_\alpha)
$$  
$$
(\nabla_{JV_\beta}\Omega)(JV_\beta,X)=\sum_{\alpha \in R_{\Theta}^+} a_\alpha (\nabla_{JV_\beta}\Omega)(JV_\beta,V_\alpha)+\sum_{\alpha \in R_{\Theta}^+} b_\alpha (\nabla_{JV_\beta}\Omega)(JV_\beta,JV_\alpha).
$$
%and
%$$
%(\delta\Omega)(X)=\sum_{\beta \in R_{\Theta}^+}\left\lbrace (\nabla_{V_\beta}\Omega)(V_\beta,X)+(\nabla_{JV_\beta}\Omega)(JV_\beta,X)\right\rbrace ,
%$$

%Now, $(\nabla_X\Omega)(Y,Z)=g(Y,(\nabla_XJ)Z)$, for $X,Y,Z \in \mf{m}$, then

Now it is sufficient to see that each term in the four sum above is zero.  In fact, since $(\nabla_X\Omega)(Y,Z)=g(Y,(\nabla_XJ)Z)$, for $X,Y,Z \in \mf{m}$, then

\begin{align*}
	(\nabla_{V_\beta}\Omega)(V_\beta,V_\alpha)&= g(V_\beta, (\nabla_{V_\beta}J)V_\alpha)\\
	&= g(V_\beta, \nabla_{V_\beta}JV_\alpha-J(\nabla_{V_\beta}V_\alpha))\\
	&=g\left(\dfrac{1}{\sqrt{2\lambda_\beta}}A_\beta , \dfrac{i}{2\sqrt{\lambda_\beta\lambda_\alpha}}\nabla_{A_\beta}S_\alpha-\dfrac{1}{2\sqrt{\lambda_\beta\lambda_\alpha}}J(\nabla_{A_\beta}A_\alpha)\right).\\
%	&=\dfrac{1}{2\sqrt{\lambda_\beta\lambda_\alpha}}\dfrac{1}{\sqrt{2\lambda_\beta}}g(A_\beta, i\nabla_{A_\beta}S_\alpha-J(\nabla_{A_\beta}A_\alpha))\\
%	&=\dfrac{1}{2\sqrt{\lambda_\beta\lambda_\alpha}}\dfrac{1}{\sqrt{2\lambda_\beta}}g(A_\beta, i(r_{\beta,\alpha}A_{\beta+\alpha}+r_{\beta,-\alpha}S_{\beta-\alpha})-r_{\beta,\alpha}J(S_{\beta+\alpha})+r_{\beta,-\alpha}J(A_{\beta-\alpha}))\\
%	&=\dfrac{1}{2\sqrt{\lambda_\beta\lambda_\alpha}}\dfrac{1}{\sqrt{2\lambda_\beta}}g(A_\beta, i(r_{\beta,\alpha}A_{\beta+\alpha}+r_{\beta,-\alpha}S_{\beta-\alpha})-r_{\beta,\alpha}i\varepsilon_{\beta+\alpha}A_{\beta+\alpha}+r_{\beta,-\alpha}i\varepsilon_{\beta-\alpha}S_{\beta-\alpha})\\
%	&=\dfrac{1}{2\sqrt{\lambda_\beta\lambda_\alpha}}\dfrac{1}{\sqrt{2\lambda_\beta}}g(A_\beta, (1-\varepsilon_{\beta+\alpha})i r_{\beta,\alpha}A_{\beta+\alpha}+(1-\varepsilon_{\beta-\alpha})i r_{\beta,-\alpha}S_{\beta-\alpha})\\
\end{align*}
From Lemma \ref{connections AS} the right side in last equation above is 
\begin{equation*}
 \dfrac{1}{2\lambda_\beta\sqrt{2\lambda_\alpha}}\left\lbrace    (1-\varepsilon_{\beta+\alpha})i r_{\beta,\alpha} g(A_\beta, A_{\beta+\alpha})+ (1-\varepsilon_{\beta-\alpha})i r_{\beta,-\alpha}g(A_\beta, S_{\beta-\alpha})\right\rbrace=0,
\end{equation*}
because $g(A_\beta, A_{\beta+\alpha})=0=g(A_\beta, S_{\beta-\alpha})$.

Analogously, 
\begin{align*}
	(\nabla_{V_\beta}\Omega)(V_\beta,JV_\alpha)	%g(V_\beta,(\nabla_{V_\beta}J)JV_\alpha)\\
	%&=g(V_\beta,\nabla_{V_\beta}JV_\alpha-J(\nabla_{V_\beta}JV_\alpha))\\
	&=\dfrac{i\varepsilon_\alpha}{2\lambda_\beta\sqrt{2\lambda_\alpha}}g(A_\beta,\nabla_{A_\beta}S_\alpha - J(\nabla_{A_\beta}S_\alpha))\\
	&=0
\end{align*}
since $g(A_\beta,\nabla_{A_\beta}S_\alpha)=g(A_\beta, J(\nabla_{A_\beta}S_\alpha))=0$.
Using Lemma \ref{connections AS}, by a similar computation, one obtains:
$$
(\nabla_{JV_\beta}\Omega)(JV_\beta,V_\alpha)=0=	(\nabla_{JV_\beta}\Omega)(JV_\beta,JV_\alpha).
$$
	
\end{proof}

%Considering the Weyl basis $\{X_\alpha,\quad \alpha \in R_\Theta\}$ of $\mf{m}^\CC$ we have
%
%\begin{proposition}
%\begin{eqnarray*}
%2(\nabla_{X_\alpha}\Omega)(X_\beta,X_\gamma)=\left\{\begin{array}
%	[c]{cc}%
%	in_{\gamma,\alpha}(\varepsilon_\gamma+\varepsilon_\beta)(\lambda_{\beta}+\lambda_\alpha-\lambda_\gamma), & \mbox{if}\quad \alpha+\beta+\gamma=0 \\ \\
%	\hspace*{-1cm}0,& \mbox{otherwise.}
%\end{array}
%\right.
%\end{eqnarray*}
%
%\end{proposition}  

%\begin{lemma}
%	For any $\alpha, \beta \in R_\Theta^+$, one has 
%	$$ (\nabla_{V_\beta}\Omega)(V_\beta,V_\alpha)=(\nabla_{V_\beta}\Omega)(V_\beta,JV_\alpha)=(\nabla_{JV_\beta}\Omega)(JV_\beta,V_\alpha)=(\nabla_{JV_\beta}\Omega)(JV_\beta,JV_\alpha)=0
%	$$	
%	In particular, for any $X\in \mf{m}$:
%	$$ (\nabla_{V_\beta}\Omega)(V_\beta,X)=(\nabla_{JV_\beta}\Omega)(JV_\beta,X)=0.
%	$$
%\end{lemma}

%As an immediate consequence we obtain

\begin{corollary}\label{flags ASK}
Any invariant almost Hermitian structure on $\F$ belongs to the almost semi K\"ahler class (ASK).
\end{corollary}

\section{Invariant almost Hermitian submanifolds of $\F$}\label{minimal iff ASK}

Let $\F=U/K_\Theta$ be a flag manifold. Consider $L \subset U$ a 2r-dimensional closed (and then compact) subgroup of $U$. Then $L\cap K_\Theta$ is a closed subgroup of $L$ and the coset space $L/L\cap K_\Theta$ is a submanifold of $\F$ such that the following diagram of principal bundles is commutative (see \cite{Hel}, p. 125)

\begin{center}
	\begin{tikzcd}
		L \arrow[r,"i", hook] \arrow[d, "\pi_{1}"]
		& U \arrow[d, "\pi"] \\
		L/L\cap K_\Theta \arrow[r,"I", hook]
		& U/K_\Theta.
	\end{tikzcd}
\end{center} 
where $\pi_1$ and $\pi$ are the natural mappings of $L$ onto $L/L\cap K_\Theta$ and $\pi$ of $U$ onto $U/K_\Theta$, respectively. The inclusion mapping of $L$ into $U$ is denoted by $i$ and $I\colon L/Q\longrightarrow U/K_\Theta$ denotes the mapping $l(L \cap U_\Theta)\mapsto i(l)K_\Theta$. We denote by $\mf{l}$  the Lie algebra of $L$. Since $\mf{l}$ is a subalgebra of the compact Lie algebra $\mf{u}$, the Killing form $B(\cdot, \cdot)$ is negative definite on $\mf{l}$. From  the decomposition $\mf{u}=\mf{k}_{\Theta}\oplus\mf{m}$ we obtain
$$
\mf{l}=\mf{l}\cap\mf{k}_{\Theta}\oplus\mf{n}
$$
where $\mf{n}\subset \mf{m}$ is the orthogonal complement of $\mf{l}\cap\mf{k}_{\Theta}$ in $\mf{l}$ with respect to $-B(\cdot, \cdot)$. 

The $U$-invariance of $g(,)$ implies $g([X,Y],Z)=g(X,[Y,Z]])$, for any $X,Y, Z$ in $\mf{u}$. In particular, for $X,Y\in \mf{l}\cap\mf{k}_{\Theta} $ and $Z\in \mf{n} $, we have $0=g([X,Y],Z)=g(X,[Y,Z]])$, then  $\operatorname{ad}(\mf{l}\cap\mf{k}_{\Theta})\mf{n}\subset\mf{n}$. The real tangent $T_o(L/L\cap K_\Theta) $ at the origin $o=L\cap K_\Theta $ is isomorphic to the real subspace $\mf{n}$.
\begin{lemma}
Let $(\F= U/K_\Theta, J,g)$ be an invariant almost Hermitian flag manifold and $L/L\cap K_\Theta$ as above, then $J(\mf{n})=\mf{n}$.
\end{lemma}
\begin{proof}
	
		Consider the following subset of the complementary roots

	\begin{equation}\label{Rn}
		R_{\mf{n}}=\left\lbrace \alpha\in R_{\Theta}^+: g(X,V_\alpha)\neq 0,\quad \text{or}\quad g(X,JV_\alpha)\neq 0  \quad \text{for some} \quad X \in \mf{n} \right\rbrace  
	\end{equation}
where the vectors $V_\alpha$ and $JV_\alpha$ were defined in (\ref{J-frame}).
	
	Then the following subset of orthonormal J-frame $\left\lbrace V_\alpha, J(V_\alpha); \quad \alpha \in R_{\mf{n}} \right\rbrace $ span the tangent space $\mf{n}$, because if $X\in \mf{n}$, 
	$$
	X=\sum_{\alpha \in R_{\mf{n}}} a_\alpha V_\alpha +\sum_{\alpha \in R_{\mf{n}}} b_\alpha JV_\alpha
	$$ 
	where $a_\alpha= g(X,V_\alpha)$ and $b_\alpha=g(X,JV_\alpha)$ are the real coefficients of this linear combination. 
   
   Thus, 
   	$$
   JX=\sum_{\alpha \in R_{\mf{n}}} a_\alpha JV_\alpha -\sum_{\alpha \in R_{\mf{n}}} b_\alpha V_\alpha 
   $$ 
   with $\alpha\in R_{\mf{n}}$. Now, for $\alpha\in R_{\Theta}^+$ and $L \in \mf{k}_{\Theta}$, $$
   g(V_\alpha,L)=0=g(JV_\alpha, L)
   $$
   which is valid, in particular, for $\alpha\in 	R_{\mf{n}}\subset R_{\Theta}^+$ and $L \in \mf{l}\cap\mf{k}_{\Theta}$.  Then $g(JX,L)=0$ for any $L\in \mf{l}\cap\mf{k}_{\Theta} $ and $J(\mf{n})\subset \mf{n}$. Since $X=J(-JX)$, one has $\mf{n}\subset J(\mf{n})$.
   This proves that $J(\mf{n})=\mf{n}$.
\end{proof}
In the sequel we continue denoting by $J$ and $g$ their restriction to $\mf{n}$.

Let $(\F= U/K_\Theta, J,g)$ be an invariant almost Hermitian flag manifold, from the above Lemma, the restrictions of $J$ and $g$ to $\mf{n} \subset \mf{m}$ define an invariant almost Hermitian structure on $L/L\cap K_\Theta $ in a such way that $(L/L\cap K_\Theta,J,g)$ is an holomorphic submanifold of $\F$. 

Now, since the vector fields are invariants, we have the following isomorphisms:

$$
\mathfrak{X}(\F)=\mf{m}=\bar{\mathfrak{X}}(\F),\quad \mathfrak{X}(L/L\cap K_\Theta)=\mf{n}
$$
%$$
%\bar{\mathfrak{X}}(\F)=\{X\in\mf{m};g(X,Y)\neq 0,\quad \text{for some}\quad Y\in \mf{n}\}
%$$
$$
\mathfrak{X}(L/L\cap K_\Theta))^\perp=\mf{n}^\perp=\{X\in\mf{m};g(X,\mf{n})= 0\}
$$

and the equation (\ref{fields perp}) becomes
\[
\mf{m}=\mf{n}\oplus \mf{n}^\perp.
\]

\begin{lemma}
	The orthogonal complement $\mf{n}^\perp$ admits a $g$-orthonormal basis $\{V_\alpha, JV_\alpha\}$, with $\alpha \in R_\Theta-R_{\mf{n}}$.
\end{lemma}
\begin{proof}
From the definition (\ref{Rn}), if $\alpha\in R_\Theta-R_{\mf{n}}$ then $g(X,V_\alpha)=g(X,JV_\alpha)=0$, for all $X\in \mf{n}$, hence $\{V_\alpha, JV_\alpha\}\subset \mf{n}^\perp$. Note that any $Y\in  \mf{m}$ can be written as 
\begin{align*}
	Y&=Y_\mf{n}+Y_{\mf{n}^\perp}\\ 
	&=\sum_{\alpha \in R_{\mf{n}}}\left(  a_\alpha V_\alpha +b_\alpha JV_\alpha \right)+ \sum_{\alpha \in R_\Theta-R_{\mf{n}}}\left(  c_\alpha V_\alpha +d_\alpha JV_\alpha \right). 
\end{align*}
By definition of $R_{\mf{n}}$ the $g$-orthogonal projection of $Y$ on $\mf{n}$ is given by
$$
Y_\mf{n}= \sum_{\alpha \in R_{\mf{n}}}\left(  a_\alpha V_\alpha +b_\alpha JV_\alpha \right)
$$
where $a_\alpha=g(Y,V_\alpha)$ and $b_\alpha=g(Y,JV_\alpha)$, then
$$
Y_{\mf{n}^\perp}=\sum_{\alpha \in R_\Theta-R_{\mf{n}}}\left(  a_\alpha V_\alpha +b_\alpha JV_\alpha \right).
$$
Thus the subset of the J-frame $\left\lbrace V_\alpha, JV_\alpha;\quad \alpha \in R_\Theta-R_{\mf{n}} \right\rbrace$, spans $\mf{n}^\perp$. Since we have an $g$-orthonormal subset, it is an invariant basis for $\mf{n}^\perp$.

It remains to prove that $R_\Theta-R_{\mf{n}}\neq \emptyset$. Above we proved that if $Y_{\mf{n}^\perp}\neq 0$ then  $R_\Theta-R_{\mf{n}}\neq \emptyset$. Thus, if $R_\Theta=R_{\mf{n}}$ then $\mf{n}^\perp=0$ and $L/L\cap K_\Theta$ has codimension $0$. So  $R_\Theta-R_{\mf{n}}\neq \emptyset$.
\end{proof}
 \begin{remark}
 The subset $\{V_\alpha, JV_\alpha\}$, with $\alpha \in R_{\mf{n}}$ of the invariant J-frame spans the space $\mf{n}$ but may not be a basis because it is possible that some of these vectors do not belong to $\mf{n}$. We will see in next section that when $L/L\cap K_\Theta$ is also a flag manifold, this subset of the invariant J-frame becomes a basis of the tangent space $\mf{n}$.
 \end{remark}

\begin{theorem}\label{ask}
 Let $(\F, J,g)$ be an invariant almost Hermitian flag manifold (witch belong to the class ASK). Then any  invariant holomorphic submanifold $L/L\cap K_\Theta$  of  $(\F, J,g)$ belongs to the class ASK. Moreover, $L/L\cap K_\Theta$ is necessarily minimal.
\end{theorem}

\begin{proof}
Consider the invariant J-frame  $\{V_\alpha, JV_\alpha\}$,  $\alpha \in R_\Theta-R_{\mf{n}}$, on $\mf{n}^\perp$.  For $X\in \mf{n}\subset \mf{m}$, the equation (\ref{same class}) sets:  
$$
({\delta}\Omega)(X)=({\delta}^\prime\Omega)(X)+ (\bar{\bar{\delta}}\Omega)(X).
$$
Now, by Corollary \ref{flags ASK}, $({\delta}\Omega)(X)=0$ and from Theorem \ref{codifferential} $(\bar{\bar{\delta}}\Omega)(X)=0$, computed on the invariant J-frame  $\{V_\alpha, JV_\alpha\}$,  $\alpha \in R_\Theta-R_{\mf{n}}$. So $L/L\cap K_\Theta$ belongs to the class ASK.

To see that $L/L\cap K_\Theta$ is minimal, let $\left\lbrace E_1,\dots, E_r,V_{\alpha_1},\dots, V_{\alpha_s}, JE_1, \dots, JE_r, JV_{\alpha_1},\dots, JV_{\alpha_s}\right\rbrace $ be an adapted invariant ortonormal J-frame  of $\mf{m}$, such that $\left\lbrace E_1,\dots, E_r, JE_1, \dots, JE_r\right\rbrace $ is an ortonormal J-frame  of $\mf{n}$. For $X\in \mf{m}$, we can write $X=X_\mf{n}+X_{\mf{n}^\perp}$ such that $X_\mf{n}\in \mf{n}$ and $X_{\mf{n}^\perp}\in\mf{n}^\perp$. Thus using Theorem \ref{codifferential}, we have 

\begin{eqnarray*}
	0=(\delta\Omega)(X)&=&\sum_{i=1}^{r}((\nabla_{E_i}\Omega)(E_i,X)+(\nabla_{JE_i}\Omega)(JE_i,X))\\&+& \sum_{j=1}^{s}(\underbrace{(\nabla_{V_{\alpha_j}}\Omega)(V_{\alpha_j},X)}_{=0}+\underbrace{(\nabla_{JV_{\alpha_j}}\Omega)(JV_{\alpha_j},X)}_{=0}\\
	&=& \sum_{i=1}^{r}((\nabla_{E_i}\Omega)(E_i,X_\mf{n}+X_{\mf{n}^\perp})+(\nabla_{JE_i}\Omega)(JE_i,X_\mf{n}+X_{\mf{n}^\perp}))
	\\
	&=& (\delta^\prime\Omega)(X_\mf{n})+(\bar{\delta}\Omega)(X_{\mf{n}^\perp})
\end{eqnarray*}

From (\ref{minimal}) we get that $(L/L\cap K_\Theta,J,g)$ is minimal since $(\delta^\prime\Omega)(X_\mf{n})=0$.

\end{proof}
 
\begin{corollary}\label{sk}
If $(\F, J,g)$ belongs to the class SK . Then any  invariant holomorphic submanifold $L/L\cap K_\Theta$  of  $\F$ belongs to the class SK. Moreover, $L/L\cap K_\Theta$ is necessarily minimal.
\end{corollary}
\section{Almost Hermitian flag submanifolds of $\F$}

Let $\F=G/P_{\Theta}$ be a flag manifold and $L$ a Lie subgroup semisimple of $G$. Then $Q:=L\cap P_\Theta$ is a parabolic subgroup of $L$. 

Thus the flag manifold $L/Q$ is a submanifold of $\F$ such that the following diagram of principal bundles is commutative

\begin{center}
	\begin{tikzcd}
		L \arrow[r, hook] \arrow[d, "Q"]
		& G \arrow[d, "P_\Theta"] \\
		L/Q \arrow[r, hook]
		& G/P_\Theta
	\end{tikzcd}
\end{center}

In terms of Lie theory we can obtain a submanifold $L/Q$ of $\F$ as following: 
let  $\Theta^\prime \subset \Sigma$ be a non empty subset of simple roots such that $\Theta^\prime \not\subset \Theta$ and denote by $	R^\prime:=\langle\Theta^\prime\rangle\cap R$ the roots spanned by $\Theta^\prime$.  Consider the following semisimple Lie subalgebra of $\mf{g}$:

\[
\mf{l}=\mf{h}^\prime\oplus\sum_{\alpha\in R^\prime}\mf{g}_\alpha
\]
where $\mf{h}^\prime= \sum_{\alpha\in R^\prime}\mathbb{C}H_{\alpha}\subset \mf{h}$ is a Cartan subalgebra of $\mf{l}$ such that $R^\prime$ is the root system associated to $(\mf{l},\mf{h}^\prime)$. The following subalgebra

\[
\mf{q}_{\Theta^\prime}=\mf{l}\cap\mf{p}_{\Theta}=\mf{h}^\prime\oplus
\sum_{\alpha\in R^{\prime^+}} \mf{g}_\alpha\oplus 
\sum_{\alpha\in  R^{\prime^-}\cap R(\Theta)^-}\mf{g}_{\alpha}
\]
is a parabolic subalgebra in $\mf{l}$ because it contains $\mf{b}^\prime=\mf{h}^\prime\oplus
\sum_{\alpha\in R^{\prime^+}} \mf{g}_\alpha $ which is a Borel subalgebra of $\mf l$.
Let $L$ be the connected Lie group with Lie algebra $\mf{l}$ and denote by $\f =\{Ad_L(l)\mf{q}_{\Theta^\prime}: \,  l\in L \}$ the set of subalgebras of $\mf l$ conjugate to $\mf{q}_{\Theta^\prime}$. Thus $\f=L/Q_{\Theta^\prime}$, where $Q_{\Theta^\prime}=\{l\in L; \quad Ad_L(l)\mf{q}_{\Theta^\prime}=\mf{q}_{\Theta^\prime}\}$ is the normalizer of $\mf{q}_{\Theta^\prime}$ in $L$.

We can see $\f$ as a compact homogeneous space following a similar way used to obtain $\F=U/K_\Theta$. Using the Weyl basis fixed previously, we see that  
$$
\mf v=i\mf h^\prime_\RR\oplus\sum_{\alpha\in R^{\prime^+}}\mf u_\alpha
$$
is a compact real form of $\mf l$, where $\mf h^\prime_\RR$ denotes the real space vector spanned by $\{H_\alpha: \,\alpha\in R^\prime\}$. Thus $V=\exp \mf v$ is the compact real form of $L$ corresponding to $\mf v$. Restricting the action of $L$ on $\f$ we see that 
$$
\f=V/K_{\Theta^\prime}
$$
where $K_{\Theta^\prime}=V\cap Q_{\Theta^\prime}$. The Lie algebra of $K_{\Theta^\prime}$ is the fixed point of the conjugate  $\tau\colon X_{\alpha}\mapsto-X_{-\alpha}$ of $\mf{l}$ restricted to $\mf{q}_{\Theta^\prime}$, given by
$$
\mf{k}_{\Theta^\prime}=\mf v\cap \mf{q}_{\Theta^\prime}=i\mf{h}^\prime_\RR \oplus \sum_{\alpha\in R^{\prime^+}\cap R(\Theta)^+}\mf{u}_\alpha.
$$

Thus the real tangent space of $\f$ at the origin is identified with
$$
T_o(V/K_{\Theta^\prime})=\mf n=\sum_{\alpha\in R^{\prime^+}\setminus R(\Theta)^+}\mf{u}_\alpha,
$$
and its complexification is 
$$
\mf{n}^\CC=\mf n \otimes \CC = \sum_{\alpha\in R^\prime\setminus R(\Theta)}\CC X_\alpha.
$$
Note that $R^\prime\setminus R(\Theta)\neq \emptyset$, since $\Theta^\prime \not\subset \Theta$. Of course $\mf{n}^\CC$ is a vector subspace of $\mf{m}^\CC$ and 
$$
\mf{m}^\CC-\mf{n}^\CC=\sum_{\alpha\in R_\Theta\setminus R^\prime}\CC X_{\alpha}.
$$
Then $\f$ is a submanifold of $\F$ of codimension $|R_\Theta|-|R^\prime\cap R_\Theta |=|R|-|R(\Theta)|-|R^\prime\cap R_\Theta |.$

From the previous section, it is clear that any invariant almost Hermitian structure $(g,J)$ on $\F$ can be restrict to $\f$ in a such way that $\f$ becomes an invariant holomorphic submanifold of $\F$. We will keep the same notation for $g$ and $J$ on $\f$.

Next, we denote by $Z_\mf{m}$ (resp. $Z_\mf{n}$) the component of $Z\in \mf{g}$ in the subspace $\mf{m}$ (resp. $\mf{n}$).
\begin{proposition}\label{comp equal}
	If $X,Y\in\mf{n} $ then $\left[X,Y \right]_\mf{n}=\left[X,Y \right]_\mf{m}$.
\end{proposition}
\begin{proof}
	Since $\left[X,Y \right]_\mf{n}=\left[X,Y \right]_\mf{m}\iff \left[X,Y \right]_\mf{n^\CC}=\left[X,Y \right]_\mf{m^\CC}$ we just need to prove that $\left[X_\alpha,X_\beta \right]_\mf{n^\CC}=\left[X_\alpha,X_\beta \right]_\mf{m^\CC}$, for $\alpha, \beta \in R^\prime\setminus R(\Theta) $.
	
	If $\alpha, \beta \in R^\prime\setminus R(\Theta) $ then $\alpha+\beta \in R$ or $\alpha+\beta \notin R$. If  $\alpha+\beta \notin R$ then $\left[X_\alpha,X_\beta \right]=0$, in particular $\left[X_\alpha,X_\beta \right]_\mf{n^\CC}=\left[X_\alpha,X_\beta \right]_\mf{m^\CC}$. If $\alpha+\beta \in R$ then $\alpha+\beta \in R^\prime$, because $R^\prime$ is a root system. Now $\alpha+\beta \in R^\prime\cap R(\Theta)$  or $\alpha+\beta \in R^\prime\setminus R(\Theta)$. In the first case  $\left[X_\alpha,X_\beta \right]_\mf{n^\CC}=0=\left[X_\alpha,X_\beta \right]_\mf{m^\CC}$. In the order case, $\left[X_\alpha,X_\beta \right]_\mf{n^\CC}=n_{\alpha,\beta}X_{\alpha+\beta}=\left[X_\alpha,X_\beta \right]_\mf{m^\CC}$.  
\end{proof}

\begin{example}
	Consider  $\mf g=\mathfrak{sl}(8,\CC)$. A Cartan subalgebra of  $\mathfrak{sl}(8,\mathbb{C})$ has the form
	$$
	\mathfrak{h}=\displaystyle\left\{\operatorname{diag}(\varepsilon_1, \ldots, \varepsilon_8): \, \varepsilon_i\in \mathbb{C} \mbox{ and } \sum_{i=1}^{8}\varepsilon_i=0\right\}.
	$$
	The root system of $\left(  \mathfrak{sl}(8,\mathbb{C}),
	\mathfrak{h}\right)  $ is given by linear maps $\alpha_{i,j}\colon \operatorname{diag}\left( \varepsilon_{1},\varepsilon
	_{2},\ldots,\varepsilon_{8}\right)
	\mapsto\varepsilon_{i}-\varepsilon_{j}$. The root system is 
	$$
	R=\left\{  \pm  \alpha_{i,j}
	,1\leq i<j\leq 8\right\}  
	$$
	and a choice of positive roots and simple roots associated are given by
	\[
	R^{+}=\left\{  \alpha_{i,j},1\leq i<j\leq
	8\right\},  \quad \Sigma= \left\{  \alpha_{i,i+1},1\leq i\leq
	8\right\}. 
	\]
	Now let $\Theta=\left\{\alpha_{1,2}, \alpha_{2,3},\alpha_{5,6}\right\}$ and $\Theta^\prime=\left\{\alpha_{1,2}, \alpha_{2,3},\alpha_{3,4},\alpha_{7,8}\right\}$ be subsets of $\Sigma$. Note that $\Theta^\prime \not\subset \Theta$ and the Dynkin diagram of $\Theta^\prime$ is not connected.
	
	We have
	
	$$
	R(\Theta):=\langle\Theta\rangle\cap R=\{\pm\alpha_{1,2}, \pm\alpha_{2,3},\pm\alpha_{1,3},\pm\alpha_{5,6}\}
	$$
	\[
	\p_{\Theta}:=\mf{b}^+\oplus\mf{g}_{\alpha_{2,1}}\oplus\mf{g}_{\alpha_{3,2}}\oplus\mf{g}_{\alpha_{3,1}}\oplus\mf{g}_{\alpha_{6,5}},
	\]
	
	where $\mf{b}^+=\mf{h}\oplus\sum\limits_{\alpha\in R^+}\mf{g}_{\alpha}$ denotes the Borel subalgebra.
	
	The real form compact of $\p_{\Theta}$ is isomorphic to
	\[
	\mf{k}_{\Theta}=\mf{su}(8)\cap\p_{\Theta}=\mf{s}(\mf{u}(3)\times \mf{u}(2)\times \mf{u}(1)^3)
	\]
	and $\F=SU(8)/S(U(3)\times U(2)\times U(1)\times U(1)\times (1))$. Now
	$$
	R^\prime:=\langle\Theta^\prime\rangle\cap R=\{\pm\alpha_{1,2}, \pm\alpha_{2,3},\pm\alpha_{1,3},\pm\alpha_{3,4},\pm\alpha_{1,4},\pm\alpha_{2,4},\pm\alpha_{7,8}\}.
	$$
	A basis for a Cartan subalgebra associated to root system $R^\prime$ is given by the matrix of the form
	\[
	H_{\alpha_{i,j}}= E_{ii}-E_{jj}\quad 1\leq i\neq j\leq 4\quad \text{and}\quad  7\leq i\neq j\leq 8
	\]
	where $E_{rs}$ denotes the matrix $8\times 8$  with $1$ at position $(r,s)$ and zeros everywhere else.
	Thus
	\[
	\mf{l}=\mf{h}^\prime\oplus\sum_{\alpha\in R^\prime}\mf{g}_\alpha=\mf{sl}(4,\CC)\times \mf{sl}(2,\CC)
	\]
	and the real compact form of $\mf{l}$ is 
	
	$$
	\mf v=i\mf h^\prime_\RR\oplus\sum_{\alpha\in R^{\prime^+}}\mf u_\alpha= \mf{s}(\mf{u}(4))\times \mf{s}(\mf{u}(2)).
	$$
	The following is a parabolic subalgebra of $\mf{l}$
	\[
	\mf{q}_{\Theta^\prime}=\mf{l}\cap\mf{p}_{\Theta}=\mf{h}^\prime\oplus
	\sum_{\alpha\in R^{\prime^+}} \mf{g}_\alpha\oplus 
	\sum_{\alpha\in  R^{\prime^-}\cap R(\Theta)^-}\mf{g}_{\alpha}=\mf{h}^\prime\oplus\sum_{\alpha\in R^{\prime^+}} \mf{g}_\alpha\oplus
	\mf{g}_{\alpha_{2,1}}\oplus\mf{g}_{\alpha_{3,2}}\oplus\mf{g}_{\alpha_{3,1}}
	\]
	its real compact form is given by
	\[
	\mf{k}_{\Theta^\prime}=\mf v\cap \mf{q}_{\Theta^\prime}=i\mf{h}^\prime_\RR \oplus \sum_{\alpha\in R^{\prime^+}\cap R(\Theta)^+}\mf{u}_\alpha= \mf{s}(\mf{u}(3)\times \mf{u}(1))\times \mf{s}(\mf{u}(1)\times \mf{u}(1)).
	\]
	
	Then $\f=S(U(4))\times S(U(2))/S(U(3)\times U(1))\times S(U(1)\times U(1))$. Thus we obtain the following commutative diagram
	
	\begin{center}
		\begin{tikzcd}
			S(U(4))\times S(U(2))\arrow[r, hook] \arrow[d]
			& SU(8) \arrow[d] \\ 
			\dfrac{S(U(4))\times S(U(2))}{S(U(3)\times U(1))\times S(U(1)\times U(1))} \arrow[r, hook]
			&\dfrac{SU(8)}{S(U(3)\times U(2)\times U(1)^3)}.
		\end{tikzcd}
	\end{center} 
\end{example}

\begin{theorem}\label{tg}
	If $(\F,g,J)$ belongs to the class ASK or SK, then any holomorphic flag submanifold $\f$ of $\F$ belongs to the same class. Moreover, $\f$ is a totally geodesic submanifold of $\F$ with respect to any invariant metric $g$ on $\F$.
	
\end{theorem}

\begin{proof}
	The first part is follows immediately from Theorem \ref{ask} and Corollary \ref{sk}.
	%an $g$-orthonormal basis on $\mf{n}$ is given by $\{V_\alpha,JV_\alpha\}$, with $\alpha \in R^{\prime^+}- R(\Theta)^+ $, then from Lemma $\ref{codifferential}$ we get $(\delta^\prime\Omega)(X)=0$ for any $X\in \mf{n}$, where $\delta^\prime$ denotes the codifferential operator on $\f$.
	
		Now, for $X,Y\in \mf{n}$, from Gauss formula, the second fundamental form is given by: 
	\begin{eqnarray*}
		\alpha(X,Y)&=&\nabla_XY-\nabla^\prime_XY\\
		&=& -\dfrac{1}{2}[X,Y]_{\mf m}+U(X,Y)+\dfrac{1}{2}[X,Y]_{\mf n}-U^\prime(X,Y).		
	\end{eqnarray*}
where $\nabla^\prime_XY$ is the tangential component of $\alpha(X,Y)$ (hence the connection on $\f$) and $U^\prime\colon \mf{n}\times \mf{n}\longrightarrow \mf{n}$ is the symmetric bilinear map defined by 

\begin{equation*}
	2g(U^\prime(X,Y),Z)= g([Z,X]_{\mf n},Y)+g(X, [Z,Y]_{\mf n}), \quad X,Y,Z\in \mf{n}.
\end{equation*}
Using the Proposition \ref{comp equal}, we see $U|_{\mf{n}\times\mf{n}}=U^\prime$ and $\alpha (X,Y)=0$.

\end{proof}

%%%%%%%%%%%%%%%%%%%%%%%

%%%%%%%%%%%%%%%%%%%%%%%
\end{document}